\documentclass{article}
\begin{document}
\title {A Necessary and Sufficient Condition for Existence of a Rational Point on an Elliptic Curve}
\author{Puyun Gao}
\maketitle
\begin{abstract}
In this paper, the proof of the existence of a rational point on an
elliptic curve is transformed into the proof of the existence of an
integer solution for a Diophantine equation. By a new formula for
calculating the number of elements in intersection of two finite
sets, a necessary and sufficient condition for existence of a
rational point on an elliptic curve is established. This condition
is different from L-function in the Birch and Swinnerton-Dyer
conjecture.
\end{abstract}
\section{Introduction}
Elliptic curves are related to number theory, geometry, cryptography
and data transmission([1-3]).An elliptic curve is given by
\begin{equation}
 y^{2}= x^{3}+ax+b
 \end{equation}
where $a$ and $b$ are two integers.
\par The Birch and Swinnerton-Dyer conjecture has not yet been proved([4]).
Birch and Swinnerton-Dyer try to prove a necessary and sufficient
condition for the existence of infinitely many rational points on an
elliptic curve.
\par The purpose of this paper is to give a necessary and sufficient condition for the
existence of a rational point on an elliptic curve. This condition
is different from the incomplete L-function in the Birch and
Swinnerton-Dyer conjecture.
\par When $a=b=0$, the elliptic
curve (1) becomes
\begin{equation}
 y^{2}=x^{3}
\end{equation}
It is clear that the set of rational points on  the elliptic curve
(2) is given by
$$
 S_{0}=\{(r^{2},r^{3})|r\in Q\}
$$
So, in this article, we always assume that $a^{2}+b^{2}$ is not
equal to zero.
\par Two formulas and a lemma given in ref. [5] are used. For convenience, we rewrite
them.
\par Let $A=\{i_{1},i_{2},\cdots,i_{K_{1}}\}$ and
$B=\{j_{1},j_{2},\cdots,j_{K_{2}}\} $ be the two finite subsets of
the set of integer numbers, where $K_{1}$ and $K_{2}$ are two
positive integers. Write
$$M_{0}=min\{i_{1},i_{2},\cdots,i_{K_{1}},j_{1},j_{2},\cdots,j_{K_{2}}\}$$
$$N_{0}=max\{i_{1},i_{2},\cdots,i_{K_{1}},j_{1},j_{2},\cdots,j_{K_{2}}\}$$
and let $N\geq N_{0}$ and $M\leq M_{0}$ be any two integers, the
number of the elements that belong to both $A$ and $B$ is give
by([5])
\begin{equation}
\chi(A,B)=\sum\limits^{K_{1}}_{k=1}\frac{(-1)^{N+i_{k}}}{(i_{k}-M)!(N-i_{k})!}\sum\limits^{N-M}_{l=0}\sum\limits^{l}_{m=0}(-1)^{m}\sigma_{m}(M,N)i^{l-m}_{k}
T_{N-M-l}
\end{equation}
where
$$\sigma_{m}(M,N)=\sum\limits_{M\leq n_{1}<n_{2}<\cdots<n_{m}\leq
N}n_{1}n_{2}\cdots n_{m}
$$
$n_{1}$, $n_{2}$, $\cdots$, $n_{m}$ are integers, and
$T_{k}=\sum\limits^{K_{2}}_{h=1}j_{h}^{k}$
\par In ref.[5], it is proved that $A\cap B$ is nonempty if and only
if $\chi(A,B)>0$. If we write
\begin{equation}
\Omega(A,B)=(-1)^{N}\sum\limits^{N-M}_{l=0}\sum\limits^{l}_{m=0}(-1)^{m}\sigma_{m}(M,N)S_{l-m}T_{N-M-l}
\end{equation}
where $S_{j}=\sum\limits^{K_{1}}_{k=1}(-1)^{i_{k}}i_{k}^{j}$, then
the following lemma holds([5]).
\par {\bf Lemma 1}.  $\chi(A,B)$ if and only
if $\Omega(A,B)>0$.
\section{The Necessary and Sufficient Condition of the
existence of rational points on elliptic curves}
\mbox{}\hspace{16pt}In this section, we shall give a necessary and
sufficient condition for the existence of a rational point on an
elliptic curve
\subsection{An equivalent proposition}
 Let $P$, $Q$, $R$ and $S$ are four integers with $P>0$ and $R>0$. If
$\left(\frac{Q}{P},\frac{S}{R}\right)$ is a rational point on the
elliptic curve (1), we have
$$\left(\frac{Q}{P},\frac{S}{R}\right)=\left(\frac{QR}{PR},\frac{SP}{PR}\right)=\left(\frac{X}{n},\frac{Y}{n}\right)$$
where $X=QR$, $Y=SP$ and $n=PR$. Substituting
$(x,y)=\left(\frac{X}{n},\frac{Y}{n}\right)$ into the Eq. (1)
yields
\begin{equation}
 nY^{2}=X^{3}+an^{2}X+bn^{3}
\end{equation}
Thus, $(X,Y)$ is an integer solution of the following Diophantine
equation
\begin{equation}
 ny^{2}=x^{3}+an^{2}x+bn^{3}
\end{equation}
\par On the other hand, if $(X,Y)$ is an integer solution of Eq.
(6), then Eq. (5) holds. By Eq. (5), we have
$$
 \left(\frac{Y}{n}\right)^{2}= \left(\frac{X}{n}\right)^{3}+a \left(\frac{X}{n}\right)+b
$$
Thus, $(\frac{X}{n},\frac{Y}{n})$ is a rational point on the
elliptic curve (1).
\par In conclusion, we get the following theorem.
\par {\bf Theorem 1.} There is a rational point on the
elliptic curve (1) if and only if there exists a positive
integer $n$ such that the Diophantine equation (6) has an
integer solution.
\par Since the elliptic curve (1) is symmetric with respect to $x$-axis, we have the following
corollary.
\par {\bf Corollary 1}. There is a rational point on the
elliptic curve (1) if and only if there exists a positive
integer $n$ such that the Diophantine equation (6) has an
integer solution on the upper plane(including the $x$ axis).
\par As we all know, the cubic equation $x^{3}+ax+b=0$ always has a real
root. If its smallest real root is $x_{0}$, then the image of the
elliptic curve (1) is on the right hand side of the line
$x=x_{0}$.
\par Based on the calculus theory, it is not difficult to prove the
following Lemma.
\par {\bf Lemma 2}. The function $f(x)=x^{3}+ax+b$ has the following
properties.
\par {\bf (1)} If $a\geq 0$, then $f(x)$ has only a unique zero
$x_{0}$ and is strictly monotone increasing and non-negative on
interval $[x_{0}, +\infty)$.
\par {\bf (2)} If $a<0$, then $f(x)$ has maximum point
$x_{-}=-\sqrt{-\frac{a}{3}}$ and minimal point
$x_{-}=\sqrt{-\frac{a}{3}}$.
\par{\bf (3)} If $a<0$ and
$b<\frac{2a}{3}\sqrt{-\frac{a}{3}}$, then $f(x)$ has only a unique
zero $x_{0}$ and is strictly monotone increasing and non-negative on
interval $[x_{0}, +\infty)$.
\par{\bf (4)} If $a<0$ and
$b=\frac{2a}{3}\sqrt{-\frac{a}{3}}$, then $f(x)$ has two different
zero $x_{0}=x_{-}=2\sqrt[3]{\frac{b}{2}}$ and
$x_{1}=-2\sqrt[3]{\frac{b}{2}}$, and is strictly monotone increasing
and non-negative on interval $[x_{1}, +\infty]$, and negative on
interval $[x_{0}, x_{1})$ except for point $x_{0}$.
\par {\bf (5)} If $a<0$ and $\frac{2a}{3}\sqrt{-\frac{a}{3}}<b<-\frac{2a}{3}\sqrt{-\frac{a}{3}}$,
then $f(x)$ has three different zero $x_{0}$, $x_{1}$ and
$x_{2}$(with $x_{0}< x_{-}<x_{1}<x_{+}<x_{2}$), and is strictly
monotone increasing and non-negative on interval $[x_{0}, x_{-}]$ or
$[x_{2}, +\infty)$, strictly monotone decreasing and non-negative on
interval $[x_{-}, x_{1}]$, and negative on interval $(x_{1},
x_{2})$.
\par {\bf (6)} If $a<0$ and $b=-\frac{2a}{3}\sqrt{-\frac{a}{3}}$, then $f(x)$
has two different zero $x_{0}$ and  $x_{1}=x_{+}$, and is strictly
monotone increasing and non-negative on interval $[x_{0}, x_{-}]$ or
$[x_{1}, +\infty)$, strictly monotone decreasing and non-negative on
interval $[x_{-}, x_{1}]$.
\par {\bf (7)} If
$a<0$ and $b>-\frac{2a}{3}\sqrt{-\frac{a}{3}}$, then $f(x)$ has only
a unique zero $x_{0}$, and is strictly monotone increasing and
non-negative on interval $[x_{0}, x_{-}]$ or $[x_{+}, +\infty)$,
strictly monotone decreasing and non-negative on interval $[x_{-},
x_{+}]$.
\par Write
$$
 f_{n}(x)=x^{3}+an^{2}x+bn^{3}
$$
which can be written as
\begin{equation}
 f_{n}(x)=n^{3}\left[\left(\frac{x}{n}\right)^{3}+a\left(\frac{x}{n}\right)+b\right]
 \end{equation}
By (7) and Lemma 2, we have the following Lemma.
\par {\bf Lemma 3}. The function $f_{n}(x)$ has the following
properties.
\par {\bf (1)} If $a\geq 0$, or $a<0$ and $b<\frac{2a}{3}\sqrt{-\frac{a}{3}}$, then $f_{n}(x)$
is strictly monotone increasing and non-negative on interval
$[nx_{0}, +\infty)$.
\par{\bf (2)} If $a<0$ and
$b=\frac{2a}{3}\sqrt{-\frac{a}{3}}$, then $f_{n}(x)$ is strictly
monotone increasing and non-negative on interval $[nx_{1},
+\infty)$, and negative on interval $[nx_{0}, nx_{1})$ except for
point $nx_{0}$.
\par {\bf (3)} If $a<0$ and
$\frac{2a}{3}\sqrt{-\frac{a}{3}}<b<-\frac{2a}{3}\sqrt{-\frac{a}{3}}$,
then $f_{n}(x)$ is strictly monotone increasing and non-negative on
interval $[nx_{0}, nx_{-}]$ or $[nx_{2}, +\infty)$, strictly
monotone decreasing and non-negative on interval $[nx_{-}, nx_{1}]$,
and negative on interval $(nx_{1}, nx_{2})$.
\par {\bf (4)} If $a<0$ and $b=-\frac{2a}{3}\sqrt{-\frac{a}{3}}$, then $f_{n}(x)$
is strictly monotone increasing and non-negative on interval
$[nx_{0}, nx_{-}]$ or $[nx_{1}, +\infty)$, strictly monotone
decreasing and non-negative on interval $[nx_{-}, nx_{1}]$.
\par {\bf (5)} If  $a<0$ and
$b>-\frac{2a}{3}\sqrt{-\frac{a}{3}}$, then $f_{n}(x)$ and is
strictly monotone increasing and non-negative on interval $[nx_{0},
nx_{-}]$ or $[nx_{+}, +\infty]$, strictly monotone decreasing and
non-negative on interval $[nx_{-}, nx_{+})$.
\par Write
$$J_{kn}=[nx_{k}]\mbox{}\hspace{16pt}k=0,1,2$$
$$J_{3n}=[nx_{+}]\mbox{}\hspace{16pt}J_{4n}=[nx_{-}]$$
where $[x]$ is the integral part of $x$.
\par Constructing the following finite
subsets of the set of integer numbers
$$A_{n}=\{ni^{2}|i=1,2,\cdots,I\}$$
$$B_{kn}=\{f_{n}(J_{kn}+j)|j=0,1,2,\cdots,J\}\mbox{}\hspace{16pt}k=0,1,2,3$$
$$B_{4n}=B_{n}([nx_{0},nx_{-}])=\{f_{n}(J_{0n}+j)|j=0,1,2,\cdots,J_{4n}-J_{0n}\}$$
$$B_{5n}=B_{n}([nx_{-},nx_{1}])=\{f_{n}(J_{4n}+j)|j=0,1,2,\cdots,J_{1n}-J_{4n}\}$$
$$B_{6n}=B_{n}([nx_{-},nx_{+}])=\{f_{n}(J_{4n}+j)|j=0,1,2,\cdots,J_{3n}-J_{4n}\}$$
where $I$ and $J$ are any two integers. Take
$$M_{kn}=min\{i^{2}n,f_{n}(J_{kn}+j)|i=1,2,\cdots,I;j=0,1,2,\cdots,J\}\mbox{}\hspace{16pt}k=0,1,2,3$$
$$N_{kn}=max\{i^{2}n,f_{n}(J_{kn}+j)|i=1,2,\cdots,I;j=0,1,2,\cdots,J\}\mbox{}\hspace{16pt}k=0,1,2,3$$
$$M_{4n}=min\{i^{2}n,f_{n}(J_{0n}+j)|i=1,2,\cdots,I;j=0,1,2,\cdots,J_{4n}-J_{0n}\}$$
$$N_{4n}=max\{i^{2}n,f_{n}(J_{0n}+j)|i=1,2,\cdots,I;j=0,1,2,\cdots,J_{4n}-J_{0n}\}$$
$$M_{5n}=min\{i^{2}n,f_{n}(J_{4n}+j)|i=1,2,\cdots,I;j=0,1,2,\cdots,J_{1n}-J_{4n}\}$$
$$N_{5n}=max\{i^{2}n,f_{n}(J_{4n}+j)|i=1,2,\cdots,I;j=0,1,2,\cdots,J_{1n}-J_{4n}\}$$
$$M_{6n}=min\{i^{2}n,f_{n}(J_{4n}+j)|i=1,2,\cdots,I;j=0,1,2,\cdots,J_{3n}-J_{4n}\}$$
$$N_{6n}=max\{i^{2}n,f_{n}(J_{4n}+j)|i=1,2,\cdots,I;j=0,1,2,\cdots,J_{3n}-J_{4n}\}$$
By (3), we have
\begin{equation}
\Omega(A_{n},B_{kn})=(-1)^{N_{kn}}\sum\limits^{N_{kn}-M_{kn}}_{l=0}\sum\limits^{l}_{m=0}(-1)^{m}\sigma_{m}(M_{kn},N_{kn})S_{l-m}^{n}T_{N-M-l}^{kn}
\end{equation}
where
$$S_{h}^{n}=n^{h}\sum\limits^{I}_{i=1}(-1)^{ni^{2}}i^{2h}$$
$$T_{l}^{kn}=\sum\limits^{J}_{j=0}\left(f_{n}(J_{kn}+j)\right)^{l}\mbox{}\hspace{16pt}k=0,1,2,3$$
$$T_{l}^{4n}=\sum\limits^{J_{4n}-J_{0n}}_{j=0}\left(f_{n}(J_{0n}+j)\right)^{l}$$
$$T_{l}^{5n}=\sum\limits^{J_{1n}-J_{4n}}_{j=0}\left(f_{n}(J_{4n}+j)\right)^{l}$$
$$T_{l}^{6n}=\sum\limits^{J_{3n}-J_{4n}}_{j=0}\left(f_{n}(J_{4n}+j)\right)^{l}$$
\par It is clear that $(0,0))$ is an integer solution of the Diophantine
equation (6) if $b=0$, and not its solution if $b\neq 0$. Thus,
the following theorem holds.
\par{\bf Theorem 2}. The sufficient and necessary condition that the
Diophantine equation (6) has an integer solution is that there
are two positive integers $I$ and $J$ such that one of the following
seven conditions holds.
\par {\bf (1)} $a\geq 0$ and $b\neq 0$, or  $a<0$ and $b<\frac{2a}{3}\sqrt{-\frac{a}{3}}$, $\Omega(A_{n},B_{0n})>0$.
\par{\bf (2)} $a<0$ and $b=\frac{2a}{3}\sqrt{-\frac{a}{3}}$, $nx_{0}$
is an integer. \par {\bf (3)} $a<0$ and
$b=\frac{2a}{3}\sqrt{-\frac{a}{3}}$, $nx_{0}$ is not an integer,
$\Omega(A_{n},B_{1n})>0$.
\par {\bf (4)} $a<0$, $b\neq 0$ and
$\frac{2a}{3}\sqrt{-\frac{a}{3}}<b<-\frac{2a}{3}\sqrt{-\frac{a}{3}}$,
$\Omega(A_{n},B_{2n})+\Omega(A_{n},B_{4n})+\Omega(A_{n},B_{5n})>0$
\par {\bf (5)} $a<0$ and $b=-\frac{2a}{3}\sqrt{-\frac{a}{3}}$, $\Omega(A_{n},B_{1n})+\Omega(A_{n},B_{4n})+\Omega(A_{n},B_{5n})>0$
\par {\bf (6)} $a<0$ and
$b>-\frac{2a}{3}\sqrt{-\frac{a}{3}}$,
$\Omega(A_{n},B_{3n})+\Omega(A_{n},B_{4n})+\Omega(A_{n},B_{6n})>0$.
\par {\bf (7)} $b=0$.
\par {\bf Corollary 2}. The sufficient and necessary condition that
the elliptic curve (1) has a rational point is that there are
three positive integers $n$, $I$ and $J$ such that one of seven
conditions in Theorem 2 holds.
\section{Conclusion}
In this paper, we give a sufficient and necessary condition for the
existence of a rational point on an elliptic curve. The condition is
a bivariate polynomial of $a$ and $b$ which is different from
incomplete L-function in the Birch and Swinnerton-Dyer conjecture.
In order to determine whether $\Omega(A_{n},B_{kn})$ is greater than
zero, the key is simplification of the formula (8). In addition,
can the formula (8) be simplified by Newton$^{,}s$ formula([6])?
\par Whether the existence of a rational point on an elliptic curve
can be proved by countet-proof. That is, assuming that
$\Omega(A_{n},B_{kn})=0$ for any three positive integers $n$, $I$
and $J$ , can we draw a contradictory conclusion?

\par College of Aerospace Science and Engineering, National University of Defense Technology,
Changsha, Hunan 410073, People$^{,}$s Republic of China
\par {\bf E-mail address: gfkdgpy@hotmail.com}
\end{document}